\documentclass[12pt]{article}

\usepackage{amssymb}
\usepackage{amsmath}
\usepackage{graphicx}
\usepackage{enumerate}
\usepackage{mathrsfs}
\usepackage{latexsym}
\usepackage{psfrag}
\usepackage{mathrsfs}
\usepackage{hyperref}

\newcommand{\qed}{\hfill $\square$\\
\vspace{0.1cm}}

\newcommand{\hb}{\bar{h}}
\newcommand{\A}{\mathcal{A}}
\newcommand{\B}{\mathcal{B}}

\newcommand{\F}{\mathcal{F}}

\newtheorem{theorem}{Theorem}[section]
\newtheorem{lemma}[theorem]{Lemma}
\newtheorem{proposition}[theorem]{Proposition}
\newtheorem{corollary}[theorem]{Corollary}

\newtheorem{example}[theorem]{Example}

\newenvironment{proof}{\noindent{\em Proof.}}{\qed}

\newcommand{\nchn}{\binom{n}{\lfloor \frac{n}{2}\rfloor}}

\newcommand{\comments}[1]{}

\begin{document}
\title{Minimum Weight Flat Antichains of Subsets}

\author{Jerrold R. Griggs
\thanks{Department of Mathematics, University of South Carolina,
Columbia, SC, USA 29208 (griggs@math.sc.edu).
Research supported in part by a grant from the Simons Foundation (\#282896 to Jerrold Griggs)
 and by a long-term visiting position at the IMA, University of Minnesota. } \and
Sven Hartmann
\thanks{Institut f{\"u}r Informatik, Technische Universit{\"a}t Clausthal, Julius--Albert--Str.~4,
38678 Clausthal--Zellerfeld, Germany (sven.hartmann@tu-clausthal.de)}   \and
Thomas Kalinowski
\thanks{School of Science and Technology, University of New England, Armidale, NSW 2351, Australia
(tkalinow@une.edu.au)}   \and
Uwe Leck
\thanks{Europa-Universit{\"a}t Flensburg, Auf dem Campus 1, 24943 Flensburg, Germany
(uwe.leck@uni-flensburg.de)}    \and
Ian T.~Roberts
\thanks{College of Education, Charles Darwin University, Darwin 0909, Australia
(ian.roberts@cdu.edu.au)}}

\date{\today}

\maketitle

\begin{abstract}
Building on classical theorems of Sperner and Kruskal-Katona,
we investigate antichains $\F$ in the Boolean lattice $B_n$ of
all subsets of $[n]:=\{1,2,\dots,n\}$, where $\F$ is flat,
meaning that it contains sets of at most two consecutive sizes,
say $\F=\mathcal{A}\cup\mathcal{B}$, where $\mathcal{A}$
contains only $k$-subsets, while $\mathcal{B}$ contains only $(k-1)$-subsets.
Moreover, we assume $\mathcal{A}$ consists of the first
$m$ $k$-subsets in squashed (colexicographic) order, while
$\mathcal{B}$ consists of all $(k-1)$-subsets not contained in the
subsets in $\mathcal{A}$.
Given reals $\alpha,\beta>0$, we say the weight of $\F$ is
$\alpha\cdot|\mathcal{A}|+\beta\cdot|\mathcal{B}|$.
We characterize the minimum weight antichains $\F$ for any
given $n,k,\alpha,\beta$, and we do the same when in addition
$\F$ is a maximal antichain.
We can then derive asymptotic results on both the minimum
size and the minimum Lubell function.

%Keywords:  Spanning trees, connected domination, Hamming code
\vfill%
{\noindent
\small
\textbf{Keywords:}\/ antichain,
                   Sperner family,
                   flat antichain,
                   Kruskal-Katona Theorem,
                   Lubell function
\\
\textbf{MSC (2010):}\/ 05D05 (primary), 06A07
\\
\textbf{Proposed Running Head:}\/ Flat Antichains
}

\end{abstract}

%\baselineskip=2.0\baselineskip

\section{Introduction}\label{sec:Over}

In the Boolean lattice $B_n$ of all subsets of $[n]:=\{1,2,\dots,n\}$ an {\em antichain}\/
is a family of subsets such that no one contains any other one.
A classical theorem of Sperner~\cite{sperner} gives the maximum size of an antichain, $\nchn$.
Let us denote by $\binom{[n]}{k}$ the collection of all $k$-subsets of $[n]$.
Sperner's Theorem also says that the only maximum-sized antichains are
$\binom{[n]}{\lfloor \frac{n}{2}\rfloor}$ and
$\binom{[n]}{\lceil \frac{n}{2}\rceil}$ (which are the same for even $n$).
Trivially, by taking $a$ subsets of size $\lfloor \frac{n}{2}\rfloor$, we obtain antichains
$\A$ of every size $a$ up to $\nchn$, and thus these are all possible antichain sizes in $B_n$.

What if we also require that $\A$ be {\em maximal}, which means that $\A$ is no longer an
antichain when any other subset is added to it?
What are the possible sizes of maximal antichains in $B_n$?
For collections of a single size $k$, we must take all of them to get a maximal antichain, and
this gives us sizes $\binom{n}{k}$.
What other sizes are possible?
Can all possible sizes be achieved by antichains with subsets of at most two different
consecutive sizes?
An antichain in $B_n$ such that all sets it contains have size $k-1$ or $k$, for some $k$,
is said to be {\em flat}.
It appears that most possible sizes of maximal antichains can indeed be achieved by flat
maximal antichains, but this is a large problem that we expect to say more about in a future
series of papers.

Another classical result, known as the Kruskal-Katona Theorem~\cite{kru,kat}, answers the following
question:  Given numbers $n,k,m$, how can one select a family $\F$ of
$m$ $k$-subsets of $[n]$ in order to minimize the size of the family $\Delta(\F)$ of
$(k-1)$-subsets that are each contained in some set in $\F$?
The family $\Delta(\F)$ is called the {\em shadow}\/ of $\F$.
We observe that the $(k-1)$-subsets {\em not}\/ in $\Delta(\F)$ form an antichain with $\F$, that is
for any family $\F\subseteq \binom{[n]}{k}$, we get a flat antichain
$\F\cup (\binom{[n]}{k-1} - \Delta(\F))$.
Note that such an antichain need not be maximal, as there could be additional $k$-subsets
of $[n]$ that contain none of the sets in $\binom{[n]}{k-1} - \Delta(\F)$.
In any case, selecting $\F$ with minimum shadow, we obtain an antichain of this type
that has maximum size for given $n,k,m$.
The study of flat antichains can be viewed as a generalization of the Kruskal-Katona
Theorem.

Another interesting measure of a family of subsets $\F\subseteq 2^{[n]}$ is called the
{\em volume}, given by $V(\F):= \sum_{F\in\F} |F|$.
We say that two antichains are {\em equivalent}\/ if they have the same size and the same volume.
This notion induces an equivalence relation on the class of all antichains in $B_n$.
It is remarkable that every antichain in $B_n$ is equivalent to some flat antichain!
The results of Kisv\"olcsey \cite{kis} and Lieby \cite{lie2,lie3} perfectly complement
each other to give the following theorem:

\begin{theorem} \label{FAT} (Flat Antichain Theorem)
Let $\F\subseteq B_n$ be an antichain. Then there is a flat antichain $\F'\subseteq B_n$
with $|\F'|=|\F|$ and $V(\F')=V(\F)$.
\end{theorem}

This easily-stated result is surprisingly hard to prove.
It says that each equivalence class of antichains contains a flat one.

In Section~2 we shall observe that the flat antichains are the extremal representatives of
their equivalence classes in a more general context where we apply a {\em weight function}\/
$w:B_n\rightarrow [0,\infty)$ to the subsets.
Moreover, we assume that there are real numbers $w_i\ge0$, $0\le i\le n$, such that
the weight of every subset $F\in B_n$ is $w_{|F|}$, so depends only on its size.
For a family $\F\subseteq 2^{[n]}$, its {\em weight}\/ $w(\F):=\sum_{F\in\F} w(F)$.
In Section~2
we shall compare the weight of a flat antichain to others in its class when there are
conditions on the convexity of $w$.

Besides the two examples of weight functions we have seen,
which are the size and the volume of a family of subsets $\F$,
there is a third measure of some interest.
Following~\cite{GriLiLu} we define the {\em Lubell function\/}
of $\F$ by $\hb(\F) := \sum_{F\in\F} \frac 1{\binom n{|F|}}$.
%%%Note original draft had extra factor of |F| in each term of sum
%%%I believe this was an error, and removed it
The notation is motivated by the fact that $\hb(\F)$ equals the average number of
times a random maximal chain meets $\F$,
as compared to the height of $\F$,
which is the maximum number of times over all full chains.
The name refers to Lubell~\cite{Lub},
who gave a very short proof of Sperner's Theorem by showing in this
notation that $\hb(\F)\le1$ for an antichain $\F$.
This is sometimes called the BLYM inequality,
which also gives credit to others who discovered a version of it~\cite{Bol,Yam,Mes}.
In recent years the Lubell function has played a
role in the study of families of subsets $\F$ that do not contain a given subposet
~\cite{GriLiLu} (see survey~\cite{GriLi}).
In Section~2 we observe that flat antichains have minimum Lubell values in their
equivalence classes.

%%%%%%%%%%%%%%%%%%%%%%%%%%%%%%%%%%%%%%%%%%%%%%%%%%%%%%%%%%%%%%%%%%%%%%%%%%%%%%%%%%%

We require some additional notions in the paper.
First suppose we have a family $\F\subseteq \binom{[n]}{k}$ for some $n,k$.
We have already defined the shadow $\Delta(\F)$.
We also require the dual notion:
The {\em shade}\/ $\nabla(\F)$ is the collection
of $(k+1)$-subsets of $[n]$ that contain some set in $\F$.
Now consider a flat antichain $\F=\A\cup\B$, where $\A\subseteq \binom{[n]}{k}$
and $\B\subseteq \binom{[n]}{k-1}$ for some $k$, $1\le k\le n$.
We say that $\F$ is {\em full flat}\/ provided that $\B=\binom{[n]}{k-1}-\Delta(\A)$.
It means that given $\A$, we put all sets we could into $\B$ to still have an antichain.
Moreover, if a full flat antichain $\F$ satisfies $\A=\binom{[n]}{k}-\nabla(\B)$,
we say it is {\em maximal flat}.
%Adding any set of size $k-1$ or $k$ to $\F$ would destroy the antichain property.
It is easily checked that a flat antichain is maximal flat if and only if it is
a maximal antichain.
%verify this is true!  done

For example, $\mathcal{F}=\bigl\{ \{1,2\}, \{1,3\}, \{4\} \bigr\}$ is a full flat
antichain in $B_4$ as all singletons other than $\{4\}$ are covered by the 2-sets in $\mathcal{F}$.
On the other hand, $\mathcal{F}$ is not maximal flat, since we would still have an antichain
if we added in $\{2,3\}$.

It is easy to see that for a flat antichain $\mathcal{F}$, its complement
$\overline{\mathcal{F}}:=\{[n]- F:F\in\mathcal{F}\}$ is also a flat antichain.
Further, a flat antichain $\F$ is maximal if and only if $\overline{\F}$ is maximal.
Moreover, if a flat antichain $\F$ is full, then $\overline{\F}$ is also full, only if
$\F$ is maximal.

The paper~\cite{ghklr} investigates the minimum size of a maximal flat antichain for given $n,k$.
This is a daunting problem in general, so the authors focus on the case of $k=3$.
They show that the minimum size of a maximal flat antichain consisting of
2-sets and 3-sets is $\binom{n}{2}-\lfloor (n+1)^2/8\rfloor$, and
all such antichains of minimum size are determined.
Note that such antichains are smaller than the smaller of the two rank sets,
$\binom{[n]}{2}$.
The authors actually solve the more general problem of minimising the total weight of a
maximal flat antichain with $k=3$, when there is a weight function $w$.
Results are given in~\cite{klr} for the general setting of maximal antichains of
subsets, where set sizes belong to some given set $K$.
The flat case corresponds to $K=\{k,k-1\}$.

In Section~3 we present our main result, which concerns finding full flat antichains
of minimum weight when the $k$-sets are chosen in a particular way.
The Kruskal-Katona Theorem \cite{kat,kru}, mentioned earlier, says that there is an
ordering of the sets in $\binom{[n]}{k}$ such that for any $m$ the first $m$ $k$-sets
have a shadow of smallest size among all $m$-element subsets of $\binom{[n]}k$.
Following Anderson \cite{and}, we say that for distinct sets $F,G\in \binom{[n]}{k}$,
$F$ precedes $G$ in {\em squashed order}\/ (also known as {\em colexicographic order}),
written $F<_S G$, whenever
$\max \bigl((F\cup G)- (F\cap G)\bigr)\in G$.
The first $m$ sets in squashed order have a minimum shadow.

Let $\F$ be a flat antichain of the form $\mathcal{F}=\mathcal{A}\cup\mathcal{B}$,
where $\mathcal{A}\subseteq\binom{[n]}k$ and $\mathcal{B}\subseteq\binom{[n]}{k-1}$ for
some $1\le k\le n$.
We say $\F$ is {\em squashed}\/ provided $\A$ consists of the
first $m$ elements of $\binom{[n]}k$ in squashed order for some $m$.
Moreover, if $\B$ is the complement of the shadow of $\A$, then $\F$
is a {\em full squashed flat antichain}\/ in $B_n$.
Clearly such an $\F$ is completely determined by $n,k,m$.

A full squashed flat antichain that is maximal flat is said to be a {\em maximal
squashed flat antichain}.
It is known that every full squashed flat antichain $\mathcal{A}\cup\mathcal{B}$
is contained in a unique maximal squashed flat antichain
$\mathcal{A}'\cup\mathcal{B}$ (see Proposition~\ref{FinM} below).

It is a remarkable fact that for any positive $s\le\binom{n}{\lfloor n/2\rfloor}$ there is
a full squashed flat antichain of size $s$ in $B_n$.
Full squashed flat antichains $\F$ are also important in connection with {\em ideals}\/
(also called {\em down-sets}), which are families of subsets that are closed under taking subsets.
Assume we have a weight function $w$ on $B_n$,
where each $i$-set has the same weight $w_i$, and
suppose further that $0\le w_0\le w_1\le\cdots\le w_n$.
Among all antichains of given size $s$ in $B_n$ there is a unique full squashed flat antichain
that generates an ideal of minimum weight (see Engel~\cite{eng}, Theorem 8.3.5).

Given $n,k$ we are interested here in minimizing the size of a full squashed flat
antichain $\F=\A\cup\B$, which means we seek $m$ to minimize this size.
However, the size of the shadow of the $m$ sets in $\A$ is known to be a
complicated function of $m$, which makes it challenging to solve this problem~\cite{fmrt}.

Our main result, Theorem~\ref{mainthm}, is the solution to this problem.
More generally, we determine the minimum weight of a full squashed flat antichain when
there is a weight function $w$.
Moreover, we solve the same problem for maximal squashed flat antichains.

In Section 4 we derive the results for squashed flat antichains for the specific
weights of interest to us, which are size, volume, and Lubell function value.
The paper concludes with open problems that extend the work here.

%In the next section, we solve the problem similar to the one in \cite{ghklr} for squashed
%FFA and squashed MFA for any $n$ and $k$.

%%%%%%%%%%%%%%%%%%%%%%%%%%%%%%%%%%%%%%%%%%%%%%%%%%%%%%%%%%%%%%%%%%%%%%%%%%%%%%%%%%%%%%%%%%%%%%%%%%ZZZZZZZZZZZZ

\section{Antichains with Convex Weight Functions}\label{sec:Convex}

We consider a weight function $w:B_n\rightarrow [0,\infty)$ as above, where each
$i$-set has the same weight $w_i\ge0$.
Recall that the sequence $\{w_i\}_{i=0}^n$ is {\em convex}\/ if $\{w_i-w_{i-1}\}_{i=1}^n$ is weakly increasing and
\emph{concave}\/ if $\{w_i-w_{i-1}\}_{i=1}^n$ is weakly decreasing.

\begin{proposition}\label{conc}
Let $w:B_n\rightarrow[0,\infty)$ be a weight function as above.
Furthermore, let $\mathcal{A}$ and $\mathcal{F}$ be antichains
in $B_n$, such that $\mathcal{F}$ is flat,
$|\mathcal{F}|=|\mathcal{A}|$ and $V(\mathcal{F})=V(\mathcal{A})$.
\begin{enumerate}
\item[\emph{(i)} ]
If the sequence $\{w_i\}_{i=0}^n$ is convex, then $w(\mathcal{F})\le w(\mathcal{A})$.
\item[\emph{(ii)} ]
If the sequence $\{w_i\}_{i=0}^n$ is concave, then $w(\mathcal{F})\ge w(\mathcal{A})$.
\end{enumerate}
\end{proposition}

\begin{proof}
We only prove part (i) here. The proof of (ii) is analogous.

Assume that $\{w_i\}_{i=0}^n$ is convex.
Let $\mathbf{a}=(a_0,a_1,\dots,a_n)$ be the profile vector of
$\mathcal{A}$, i.e., $a_i=|\{A\in\mathcal{A}:|A|=i\}|$.
Furthermore, let $\ell=\min\{i:a_i\ne 0\}$ and $u=\max\{i:a_i\ne 0\}$.
The weight of $\mathcal{A}$ is determined by $\mathbf{a}$, and one has
$w(\mathcal{A})=\sum_{i=\ell}^u a_iw_i$ for which we also write $w(\mathbf{a})$.
If $u-\ell\le 1$, then $\mathcal{A}$ is flat.
As $\mathcal{A}$ and $\mathcal{F}$ have the same size
and the same volume, their profile vectors must then be equal, and it follows that
$w(\mathcal{A})=w(\mathcal{F})$.
Hence, without loss of generality we can assume that $u-\ell\ge 2$.

Consider the vector $\mathbf{a}'$ obtained from $\mathbf{a}$ replacing $a_\ell$ by $a_\ell-1$,
$a_{\ell+1}$ by $a_{\ell+1}+1$, $a_u$ by $a_u-1$, and $a_{u-1}$ by $a_{u-1}+1$.
(That is, if $u-\ell=2$, then $a_{\ell+1}=a_{u-1}$ will be increased by 2.)
Note that $\sum a_i'=\sum a_i$, $\sum a_i' i=\sum a_i i$, and
\[ w(\mathbf{a})-w(\mathbf{a}')=(w_u-w_{u-1})-(w_{\ell+1}-w_\ell). \]
As $\{w_i\}$ is convex, it follows that $w(\mathbf{a}')\le w(\mathbf{a})$.
(It should be pointed out that we do not claim nor need that $\mathbf{a}'$ is the profile vector of some
antichain in $B_n$.)
Iterating this process, we transform $\mathbf{a}$ into the profile vector $\mathbf{f}$ of $\mathcal{F}$, since
$\mathcal{A}$ and $\mathcal{F}$ agree in size and volume.
This implies $w(\mathcal{F})=w(\mathbf{f})\le w(\mathbf{a})=w(\mathcal{A})$.
\end{proof}

\bigskip

Recall the Lubell function $\hb(\F)$
we mentioned in the introduction.
Proposition~\ref{conc} implies an interesting observation about flat antichains.

\begin{corollary}
Flat antichains minimize the Lubell function within their equivalence classes.
\end{corollary}

\begin{proof}
The claim follows from Proposition \ref{conc} and the straightforward fact that the sequence
$\left\{1/\binom n i\right\}_{i=0}^n$ is convex.
\end{proof}

%%%%%%%%%%%%%%%%%%%%%%%%%%%%%%%%%%%%%%%%%%%%%%%%%%ZZZZZZZZZZZZZZZZZZZZZZZZZZZZZ

% last change on Jul 17, 2019
%             by Uwe Leck

\section{Full and Maximal Squashed Flat Antichains of Minimum Weight}\label{sec:squashed}

We first review facts about the squashed (colexicographic) order and full squashed flat antichains.
For $1\le k\le n$ and $0\le m\le\binom n k$, the {\em $k$-cascade representation}\/ of $m$ is a
representation of $m$ in the form
\begin{equation}\label{cascade}
m= \sum_{i=1}^k \binom{a_i}i \text{~~with~~} a_k>a_{k-1}>\cdots >a_t\ge t>0=a_{t-1}=\cdots =a_1.
\end{equation}
The terms $\binom{a_i}i$ with $a_i=0$ could clearly be removed from the above representation
of $m$. Their only purpose here is that they will allow us a more compact formulation of
the main result (Theorem \ref{mainthm}).
It is easy to see (cf.~\cite{kat}) that for given $k$ and $m$ there is a unique $k$-cascade representation of $m$.
Moreover, if $\mathcal{A}$ is the family of the first $m$ $k$-sets in squashed order and
(\ref{cascade}) is the $k$-cascade representation of $m$, then
\begin{equation}\label{sqshad}
|\Delta\mathcal{A}| = \sum_{i=t}^k \binom{a_i}{i-1} .
\end{equation}
As noted in the introduction,
%%By the Kruskal-Katona Theorem \cite{kat,kru},
the family $\mathcal{A}$ of the first $m$ $k$-sets
in squashed order has a shadow of smallest size among all $m$-element subsets of $\binom{[n]}k$
\cite{kat,kru}.

%%Although the following should be well-known, we will provide a proof for completeness.
%%can find no reference

\begin{proposition}\label{a1=0}
Let $k\ge2$.  Let $\mathcal{F}=\mathcal{A}\cup\mathcal{B}$ be a full squashed flat antichain
with $\mathcal{A}\subseteq\binom{[n]}k$ and $\mathcal{B}\subseteq\binom{[n]}{k-1}$.
Let $m:=|\mathcal{A}|$ be represented as in {\em (\ref{cascade})}.
Then $\mathcal{F}$ is maximal flat if and only if $a_1=0$.
\end{proposition}

Our proof of Proposition \ref{a1=0} makes use of the following lemma.

\begin{lemma}\label{lem:simple_char}
Let $\mathcal F$, $\mathcal A$, $\mathcal B$, and $m$ be as in Proposition \ref{a1=0}.
Furthermore, assume that $\mathcal A = \{A_1,\ldots,A_m\}$, where the sets are listed in squashed order.
Let $A_m=\{x_1,x_2,\ldots,x_k\}$ with $x_1<x_2<\cdots<x_k$.
Then $\mathcal{F}$ is maximal flat if and only if $x_2=x_1+1$.
\end{lemma}

\begin{proof}
Put $\mathcal B = \{B_1,\ldots,B_\ell\}$, where the sets are listed in squashed order.
Then $B_1$ is the successor of $\{x_2,\ldots,x_k\}$ in $\binom{[n]}{k-1}$ with respect to squashed order.

Let $i$ be the largest index with $x_i=x_1+i-1$, i.e., $x_2=x_1+1$ if and only if $i\ge 2$.
If $x_i=n$, then $\mathcal A=\binom{[n]}{k}$ and we are done. So assume that $x_i<n$.

If $i\ge 2$, then $B_1=\{1,2,\ldots,i-2,x_{i}+1,x_{i+1},x_{i+2}\ldots,x_k\}$, and $\nabla\mathcal B$ contains
\[\{1,2,\ldots,i-2,i-1,x_{i}+1,x_{i+1},x_{i+2}\ldots,x_k\},\]
which is the successor of $A_m$ in squashed order.

If $i=1$, then the successor of $A_m$ in squashed order on $\binom{[n]}{k}$ is $\{x_1+1,x_2,\ldots,x_k\}$,
which is not contained in $\nabla\mathcal B$, because $B_1$ comes after $\{x_2,\ldots,x_k\}$.
Hence, every element of $\nabla\mathcal B$ comes after $\{x_2-1,x_2,x_3,\ldots,x_k\}$.
\end{proof}

\medskip

\noindent\textbf{Proof (Proposition \ref{a1=0}): }
The $k$-cascade representation of $m$ yields that
\begin{multline*}
\mathcal A=\binom{[a_k]}{k}\cup\left\{A\cup\{a_k+1\}\ :\ A\in\binom{[a_{k-1}]}{k-1}\right\} \\
\cup\left\{A\cup\{a_{k-1}+1,a_k+1\}\ :\ A\in\binom{[a_{k-2}]}{k-2}\right\}\\
\cup\cdots\cup\left\{A\cup\{a_2+1,a_3+1,\cdots,a_k+1\}\ :\ A\in\binom{[a_{1}]}{1}\right\}.
\end{multline*}
Let $i$ be the smallest index with $a_i>0$.
The last element of $\mathcal A$ with respect to squashed order is $A_m=\{a_i-i+1,\ldots,a_i-1,a_i,a_{i+1}+1,\ldots,a_k+1\}$.
If $i=1$, then $A_m$ starts with $a_1,\,a_{2}+1,\ldots$, and if $i>1$, then $A_m$ starts with $a_i-i+1,\,a_i-i+2$.
Now the claim follows by Lemma~\ref{lem:simple_char}.
\hfill $_\blacksquare$

\begin{proposition}\label{FinM}
Let $\mathcal{F}=\mathcal{A}\cup\mathcal{B}$ be a full squashed flat antichain
 with $\mathcal{A}\subseteq\binom{[n]}k$ and
$\mathcal{B}\subseteq\binom{[n]}{k-1}$, where $k\ge 2$.
Then $\mathcal F$ is contained in a unique maximal squashed flat antichain
$\mathcal{F}'\subseteq\binom nk\cup \binom{[n]}{k-1}$, and
$\mathcal{F}'$ is of the form $\mathcal{A}'\cup\mathcal{B}$ with
$\mathcal{A}\subseteq\mathcal{A}'\subseteq\binom{[n]}k$.
\end{proposition}

\begin{proof}
Let $\mathcal{A'}$ be the largest initial segment of $\binom{[n]}k$ with respect to squashed order such that
$\Delta\mathcal{A}'=\Delta\mathcal A$.
Clearly, $\mathcal{F}'=\mathcal{A}'\cup\mathcal{B}$ is a maximal squashed flat antichain,
and $\mathcal{A}\subseteq\mathcal{A}'$.

By the choice of $\mathcal{A}'$, any initial segment of $\binom{[n]}k$ larger than $\A'$ has its
shadow intersect $\B$, so its union with $\B$ is not an antichain.
%removed old paragraph here
\end{proof}

\bigskip

Our main result is the following characterization of all
full squashed flat antichains of minimum weight.
To avoid certain technicalities, the trivial cases $k=1$ and $k=n$ are excluded.

\begin{theorem}\label{mainthm}
Let $1< k < n$ be integers.
%Let $\alpha,\beta$ be positive real numbers, and $\lambda:=\beta/\alpha$.
Let $\mathcal{F}=\mathcal{A}\cup\mathcal{B}$ with
$\mathcal{A}\subseteq\binom{[n]}k$ and $\mathcal{B}\subseteq\binom{[n]}{k-1}$
be a full squashed flat antichain, and
let {\em(\ref{cascade})} be the $k$-cascade representation of
$m:=|\mathcal{A}|$.
Let $w$ be a weight function on $B_n$, where each $k$-set ($(k-1)$-set, resp.)
has real weight $\alpha>0$ ($\beta>0$, resp.).
Define $\lambda:=\beta/\alpha$.
Then $\mathcal{F}$ has minimum weight
$w(\mathcal{F})=\alpha\cdot |\mathcal{A}|+\beta\cdot |\mathcal{B}|$
among all full squashed flat antichains in $\binom{[n]}k\cup\binom{[n]}{k-1}$
if and only if
\[
a_i~=
\left\{\begin{array}{rl}
n-k-1+i & \text{if~~} i>1+(n-k)/\lambda , \\[.5ex]
\lceil(i-1)(\lambda+1)-1\rceil \text{~or~} \lfloor(i-1)(\lambda+1)\rfloor & \text{if~~} 1+(n-k)/\lambda
\ge i\ge 1+2/\lambda , \\[.5ex]
i & \text{if~~} 1+2/\lambda >i> 1/\lambda , \\[.5ex]
0 \text{~or~} i & \text{if~~} 1/\lambda =i , \\[.5ex]
0 & \text{if~~} 1/\lambda >i .
\end{array}\right.
\]
\end{theorem}

\begin{proof}
First, observe that with $g(m) := m-\lambda |\Delta\mathcal{A}|$ we have
\[ w(\mathcal{F}) = \alpha\cdot g(m)+\beta\binom{n}{k-1}. \]
Hence, our problem of minimizing $w(\mathcal{F})$ is equivalent to minimizing
$g(m)$ over all $m\in\left\{0,1,\dots,\binom n k \right\}$.

If $m\in\left\{\binom n k -1,\binom n k \right\}$, then $\Delta\mathcal{A}=\binom{[n]}{k-1}$ holds.
Consequently, $m=\binom n k$ does not minimize $g(m)$, and we can assume that
$m<\binom n k$, i.e.~that $a_k\le n-1$. As $a_k>a_{k-1}>\cdots >a_t$, this implies
\begin{equation}\label{aicond}
a_i=0~~\text{or}~~i\le a_i\le n-1-k+i \qquad\qquad\text{for}~~i\in [k] .
\end{equation}

By (\ref{sqshad}), we have
\begin{equation}\label{g=sumh}
g(m) = \sum_{i=t}^k h_i(a_i) ,
\end{equation}
where, for $i\in [k]$, the polynomial $h_i : \mathbb{R}\mapsto\mathbb{R}$
is defined by
\[ h_i(x)~:=~ \binom x i - \lambda\binom x{i-1}
   ~=~ \left\{\begin{array}{rl}
              x-\lambda & \text{if }~ i=1 , \\[.5ex]
              \displaystyle \frac{x+1-i(\lambda+1)}{i!} \, \prod_{j=0}^{i-2} (x-j) & \text{if }~ i\ge 2 .
              \end{array}\right. \]

Our strategy is as follows:
For each $i\in [k]$, we determine those $x\in [i,\, n-1-k+i]\cap\mathbb{Z}$ for
which $h_i(x)$ is smallest possible.
For such $x$, we choose $a_i=x$ or $a_i=0$ if $h_i(x)$ is negative or positive,
respectively. If $h_i(x)=0$, we choose $a_i\in\{0,x\}$.
Eventually, we will verify that, with the $a_i$'s chosen as described, we obtain
a proper $k$-cascade representation (\ref{cascade}), i.e., that the following
implication is true:
\begin{equation}\label{aiincr}
(i\in [k-1])~ \wedge ~(a_i>0) \quad \Longrightarrow \quad (a_i<a_{i+1}).
\end{equation}

To begin with, note that $h_1(x)=x-\lambda$ attains its global
minimum with respect to the interval $[1,n-k]$ at $x=1$, and we
have $h_1(1)=1-\lambda$.

Let $i\in [k]- \{1\}$. Then $h_i$ is a polynomial of degree $i$ with
leading coefficient $1/i!$ and zeros $0,1,\dots, i-2$ and $i(\lambda+1)-1$.
That means, $h_i(x)$ is positive and strictly increasing for $x>i(\lambda+1)-1$,
and $h_i(x)<0$ for $i-2<x<i(\lambda+1)-1$. Moreover, $h_i$ is strictly convex on
$I:=\bigl(i-2,\,i(\lambda+1)-1\bigr)$. The numbers $u:=(i-1)(\lambda+1)$
and $u-1$ both lie in $I$, and one can easily check that $h_i(u-1)=h_i(u)$.

Based on the discussion above, we distinguish three cases to find the
global minimum of $h_i(x)$ over all $x\in [i,\, n-1-k+i]\cap\mathbb{Z}$.

\emph{Case 1:~}
Assume that $u-1<i$. Note that this is equivalent to $i<1+2/\lambda$.
In this case, $h_i(x)$ is a minimum only at $x=i$, and $h_i(i)$ is positive
if $i<1/\lambda$, equals 0 if $i=1/\lambda$ and is negative if $i>1/\lambda$.

\emph{Case 2:~}
Assume that $i\le u-1$ and that $u\le n-1-k+i$. Note that this is equivalent
to $1+2/\lambda\le i\le 1+(n-k)/\lambda$.
In this case, $h_i(x)$ attains its minimum exactly for
$x\in\{ \lceil u-1 \rceil , \lfloor u \rfloor \}$, and this minimum is
negative.

\emph{Case 3:~}
Assume that $n-1-k+i<u$. Note that this is equivalent to
$1+(n-k)/\lambda <i$.
In this case, $h_i(x)$ is a minimum only at $x=n-1-k+i$ and
$h_i(n-1-k+i)<0$.

By the results of the case-by-case analysis above and (\ref{g=sumh}), $g(m)$ becomes a
minimum when the $a_i$'s are chosen as in the theorem, where the minimization
is over all choices satisfying (\ref{aicond}). Finally, a straightforward
calculation shows that (\ref{aiincr}) holds for the $a_i$'s as in the theorem.
\end{proof}

\bigskip

Note that, by Proposition \ref{a1=0}, for $\lambda<1$ the optimal
full squashed flat antichain in Theorem \ref{mainthm} is maximal.
In general, Theorem \ref{mainthm} and its proof yield the following characterization of
minimum weight maximal squashed flat antichains.

\begin{corollary}\label{minMFSA}
Let $1 < k < n$ be integers, $\alpha,\beta$ positive real numbers and
$\lambda:=\beta/\alpha$.
Furthermore, let $\mathcal{F}=\mathcal{A}\cup\mathcal{B}$ with
$\mathcal{A}\subseteq\binom{[n]}k$ and $\mathcal{B}\subseteq\binom{[n]}{k-1}$
be a maximal squashed flat antichain, and let
\emph{(\ref{cascade})} be the $k$-cascade representation of $m:=|\mathcal{A}|$.
Then $\mathcal{F}$ has minimum weight
$w(\mathcal{F})=\alpha\cdot |\mathcal{A}|+\beta\cdot |\mathcal{B}|$
among all maximal squashed flat antichains in $\binom{[n]}k\cup\binom{[n]}{k-1}$
if and only if
\begin{enumerate}
\item[\emph{(a)} ]
$\lambda\le n-k+1$ and
\[
a_i~=
\left\{\begin{array}{rl}
n-k-1+i & \text{if~~} i>1+(n-k)/\lambda , \\[.5ex]
\lceil(i-1)(\lambda+1)-1\rceil \text{~or~} \lfloor(i-1)(\lambda+1)\rfloor & \text{if~~}
1+(n-k)/\lambda \ge i\ge 1+2/\lambda , \\[.5ex]
i & \text{if~~} 1+2/\lambda >i> \max\{1/\lambda,1\} , \\[.5ex]
0 \text{~or~} i & \text{if~~} 1/\lambda =i > 1, \\[.5ex]
0 & \text{otherwise},
\end{array}\right.
\]
or
\item[\emph{(b)} ]
$\lambda\ge n-k+1$ and $a_i=0$ for $i=1,\dots,k-1$, while $a_k=n$.
\end{enumerate}
\end{corollary}

\begin{proof}
In the beginning of the proof of Theorem \ref{mainthm} we ruled out the case that $\mathcal{F}=\binom{[n]}k$ when
looking for full squashed flat antichains of minimum weight.
For $\lambda<n-k+1$, the maximal squashed flat antichain $\binom{[n]}k$ cannot be one of minimum weight
either.
This follows from the simple observation that in this case, the maximal squashed flat antichain
\[ \left(\binom{[n]}k-\nabla\bigl\{\{n-k+2, n-k+3, \dots, n\}\bigr\}\right)\cup\bigl\{\{n-k+2,n-k+3,\dots,n\}\bigr\} \]
has a smaller weight.
Now the $a_i$'s are determined as in the proof of Theorem \ref{mainthm}, with the exception
that $a_1$ must be $0$ by Proposition \ref{a1=0}. This proves the claim for $\lambda<n-k+1$.

If $\lambda>n-k+1$, then $\binom{[n]}k$ is the unique maximal squashed flat antichain of minimum weight.
To see this, assume that
$\mathcal{B}\ne\emptyset$, and use $|\nabla\mathcal{B}|\le (n-k+1)|\mathcal{B}|$, which implies that
\[ \binom{[n]}k=(\mathcal{F}-\mathcal{B})\cup\nabla\mathcal{B} \]
has a smaller weight than $\mathcal{F}$.

Finally, if $\lambda=n-k+1$, then choosing $a_1=0$ and the other $a_i$'s as in Theorem \ref{mainthm}
(i.e., $a_i=n-k-1+i$ for $i=2,\dots,k$) gives a maximal squashed flat antichain that has the
same weight as $\binom{[n]}k$.
\end{proof}

%%%%%%%%%???????????????????  do we want entire example in theorem font, not text??????
%%%%%%%%%%%%%%%%%%   Not sure, looks okay to me now
\begin{example}
Consider $n=8$ and $k=6$.
We seek all full and all maximal squashed flat antichains
in $\binom{[8]}{6}\cup\binom{[8]}{5}$ with minimum Lubell function value.
%%%That is, $n=8$ and $k=6$.
We have weights $\alpha=1/\binom{8}{6}=1/28$ and $\beta=1/\binom{8}{5}=1/56$, so that
$\lambda=\beta/\alpha=1/2$.

As $1+(n-k)/\lambda=5<6$, Theorem \ref{mainthm} yields $a_6=n-k-1+6=7$.

By $1+2/\lambda=5$, we obtain that $a_5=\lceil 4\cdot\frac 32-1\rceil=5$ or
$a_5=\lfloor 4\cdot\frac 32\rfloor=6$.

Finally, $1/\lambda=2$ implies $a_4=4$, $a_3=3$, $a_2\in\{0,2\}$, and $a_1=0$.

As we have two choices for $a_5$ and $a_2$, respectively, there are four optimal
full squashed flat antichains.
By $a_1=0$, all of them are also maximal.  We verify their Lubell values:
\begin{enumerate}
\item
With $a_5=6$ and $a_2=2$, the number of 6-sets is
\[ \textstyle |\mathcal{A}| ~=~ \binom 76 +\binom 65 + \binom 44 + \binom 33 + \binom 22 ~=~ 16, \]
while the number of 5-sets is
\[ \textstyle |\mathcal{B}| ~=~ \binom 85 - \left[ \binom 75 +\binom 64 + \binom 43 +
\binom 32 + \binom 21 \right] ~=~ 56-45 ~=~ 11 . \]
The corresponding full squashed flat antichain is the union of
$\mathcal{A}$, the collection of the first sixteen 6-sets in squashed order, and
$\mathcal{B}=\binom{[8]}{5}-\Delta\mathcal{A}$.
Its Lubell function value is
\[ \textstyle \hb(\mathcal{A}\cup\mathcal{B}) ~=~ \frac{16}{28}+\frac{11}{56} ~=~ \frac{43}{56} . \]
\item
With $a_5=6$ and $a_2=0$ we obtain
\[ \left.\begin{array}{rcl}
\textstyle |\mathcal{A}| & = & \binom 76 +\binom 65 + \binom 44 + \binom 33 ~=~ 15, \\[1ex]
\textstyle |\mathcal{B}| & = & \binom 85 - \left[ \binom 75 +\binom 64 + \binom 43 +
\binom 32 \right] ~=~ 56-43 ~=~ 13, \\[1ex]
\textstyle \hb(\mathcal{A}\cup\mathcal{B}) & = & \frac{15}{28}+\frac{13}{56} ~=~ \frac{43}{56} .
\end{array}\right. \]
\item
With $a_5=5$ and $a_2=2$ we obtain
\[ \left.\begin{array}{rcl}
\textstyle |\mathcal{A}| & = & \binom 76 +\binom 55 + \binom 44 + \binom 33 + \binom 22 ~=~ 11, \\[1ex]
\textstyle |\mathcal{B}| & = & \binom 85 - \left[ \binom 75 +\binom 54 + \binom 43 + \binom 32 +
\binom 21 \right] ~=~ 56-35 ~=~ 21, \\[1ex]
\textstyle \hb(\mathcal{A}\cup\mathcal{B}) & = & \frac{11}{28}+\frac{21}{56} ~=~ \frac{43}{56} .
\end{array}\right. \]
\item
With $a_5=5$ and $a_2=0$ we obtain
\[ \left.\begin{array}{rcl}
\textstyle |\mathcal{A}| & = & \binom 76 +\binom 55 + \binom 44 + \binom 33 ~=~ 10, \\[1ex]
\textstyle |\mathcal{B}| & = & \binom 85 - \left[ \binom 75 +\binom 54 + \binom 43 + \binom 32 \right] ~=~ 56-33 ~=~ 23, \\[1ex]
\textstyle \hb(\mathcal{A}\cup\mathcal{B}) & = & \frac{10}{28}+\frac{23}{56} ~=~ \frac{43}{56} .
\end{array}\right. \]
\end{enumerate}
\vspace*{-2ex}~\hfill
%%%%%$_\blacksquare$
\end{example}

%%%%%%%%%%%%%%%%%%%%%%%%%%%%%%%%%%%%%%%%%%%%%%%%%%%%%%%%%%%%ZZZZZZZZZZZZZZZZZ  end ffsa section
%%section cases of special interest
% last change on Jul 17, 2019
%             by Uwe Leck

\section{Cases of Special Interest}\label{sec:special}

For the smallest weight of a full squashed flat antichain and
of a maximal squashed flat antichain in $\binom{[n]}k\cup\binom{[n]}{k-1}$,
we apply
Theorem \ref{mainthm} and Corollary \ref{minMFSA}, together with (\ref{sqshad}) and
$\mathcal{B}=\binom{[n]}{k-1}-\Delta\mathcal{A}$, to obtain the formula
\begin{equation}\label{minweight}
\textstyle w(\mathcal{F})=\beta\binom n{k-1}+\sum_{i=1}^k\biggl(\alpha\binom{a_i}i-\beta\binom{a_i}{i-1}\biggr),
\end{equation}
where the $a_i$'s are chosen as in the theorem and the corollary, respectively.
(Note that for our formula to be accurate we have to adopt the somewhat unusual convention that $\binom 00$ is
0.)
We now detail what this means for minimum size, volume, and Lubell function value.

%%%%%%%%%%%%%%%%%%%%%%%%%%%%%%%%%%%%%%%%%%%%ZZZZZZZZZZZZZZZZZZZZZ
\subsection{Squashed Flat Antichains with Minimum Size}

Let $s(n,k)$ denote the minimum size of a full squashed flat antichain in $\binom{[n]}k\cup\binom{[n]}{k-1}$.
By Theorem \ref{mainthm} with $\alpha=\beta=1$, $s(n,k)$ is equal to the right-hand side of (\ref{minweight}) for
\[ a_i~=\left\{\begin{array}{rl}
        n-k-1+i & \text{if~~} i\ge n-k+2 , \\[.5ex]
        2i-3 \text{~or~} 2i-2 & \text{if~~} n-k+1\ge i\ge 3 , \\[.5ex]
        2 & \text{if~~} i=2 , \\[.5ex]
        0 \text{~or~} 1 & \text{if~~} i=1 .
        \end{array}
        \right. \]
For the particular case of a maximal squashed flat antichain we obtain minimum size
by choosing $a_1=0$ and the other $a_i$'s as above by Corollary~\ref{minMFSA}.
Consequently, the minimum size of a maximal squashed flat antichain
in $\binom{[n]}k\cup\binom{[n]}{k-1}$ is also given by the right-hand side of (\ref{minweight}).

In the display above we see that, depending on $i$, the number of possible values of $a_i$ is either 1 or 2.
We can then deduce the number of minimum size full squashed flat antichains:
For $n\le 2k-1$ there are $2^{n-k}$, while for $n\ge 2k$, there are $2^{k-1}$.
For minimum size maximal squashed flat antichains, we must have $a_1=0$, and so
in both cases $(n,k)$, half of the minimum size full squashed
flat antichains are maximal.

%%%%%%%%%%%%%%%%%%%%%%%%%%%%%%%%%%%%%%%%%%%%%%%using ian ideas
Looking at this more closely, let $\F(m)$ denote the full squashed flat antichain with
$m$ $k$-sets.
In view of Proposition~\ref{a1=0}, if $\F(m)$ is a minimum size full squashed flat antichain
with $k$-cascade parameters $a_i$, then $\F(m)$ is a maximal flat antichain when
$a_1=0$;  Otherwise, $a_1=1$ and $\F(m)$ is non-maximal.
So, the minimum size full squashed flat antichains come in pairs, $\F(m)$ and
$\F(m+1)$, such that $\F(m)$ is maximal, while $\F(m+1)$ is non-maximal.

%%%%%%%%%%%%%%%%%%%%%%%%%%%%

Using the above values for the $a_i$'s in (\ref{minweight}) gives the following formula for $s(n,k)$.

\begin{corollary}\label{sizeformula}
Let $1\le k\le (n+1)/2$. Then
\[ s(n,k)~=~s(n,n-k+1)~=~\binom n{k-1}-\sum_{i=1}^{k-1}\frac 1{i+1}\binom{2i}i . \]
\end{corollary}

Corollary \ref{sizeformula} implies that as $n$ gets larger for fixed $k$ the optimal
full squashed flat antichains look more and more like
the $(k-1)$-st level of $B_n$.
According to the next corollary, this remains true if we allow $k$ to depend on $n$.

%%%%%%%%%%%%changed restriction to n/2    and shortened proof as per Thomas
\begin{corollary}
  For any $k=k(n)\le n/2$ one has \/
  $\displaystyle s(n,k)=\bigl(1+o(1)\bigr)\binom{n}{k-1}$.
\end{corollary}

\begin{proof}
Consider
\[ S_{k-1}:=\sum_{i=1}^{k-1}\frac{1}{i+1}\binom{2i}{i},\]
which is the sum of the first $k-1$ Catalan numbers.
Topley~\cite{top} gives an upper bound on $S_{k-1}$ that is growing like a constant
times $4^k/k^{3/2}$.  This is a factor of $k$ slower than $\binom{2k}{k}$, which grows
as a constant times $4^k/k^{1/2}$, by Stirling's formula.  The claimed growth of
$s(n,k)$ follows.

%\qedhere\]

\end{proof}

%%%%%%%%%%%%%%%%%%%%%%%%%%%%%%%%%%%%%%%%%%%%%%%%%%%ZZZZZZZZZZZ
\subsection{Squashed Flat Antichains with Minimum Volume}

Using $\alpha=k$ and $\beta=k-1$ in Theorem \ref{mainthm} gives a characterization of
all minimum volume full squashed flat antichains in $\binom{[n]}k\cup\binom{[n]}{k-1}$.
In fact, working through the details, one obtains that the volume minimizer is unique:
It coincides with the first full squashed flat antichain of minimum size, and since
it has $a_1=0$, this antichain is also maximal flat.

%%%%%%%%%%%%%%%%%%%%%%%%%%%%%%%%%%%%%%%%%%%%%%%%%%%%%%%ZZZZZZZZZZZZZ
\subsection{Squashed Flat Antichains with Minimum Lubell Function Value}

To find the squashed flat antichains in $\binom{[n]}k\cup\binom{[n]}{k-1}$ that
minimize the Lubell function, we can use (\ref{minweight}) with $\alpha=1/\binom nk$ and
$\beta=1/\binom n{k-1}$, which means that $\lambda=(n-k+1)/k$.
Applying Theorem~\ref{mainthm} and simplifying the resulting expressions leads to the following results,
where the squashed flat antichains are ordered according to the number of $k$-sets in them.

\begin{description}
\item[$n\geq 2k+1$] All full (maximal, resp.) squashed flat antichains with minimum Lubell value come after the last full (maximal, resp.) squashed flat antichain with minimum size. The full squashed flat antichains with minimum Lubell value are not maximal.
\item[$n=2k$] There is a unique full (maximal, resp.) squashed flat antichain of minimum Lubell value. It is the last full (maximal, resp.) squashed flat antichain of minimum size.  Note the last maximal one is followed immediately by the last full one.
\item[$n=2k-1$] We have $\lambda=1$, and minimizing the Lubell value is equivalent to minimizing the size.  Then the full (maximal, resp.) squashed flat antichains with minimum Lubell value are those with minimum size.  There are $2^{k-1}$ full squashed flat antichains with minimum Lubell value, and half of them are maximal.  The first minimum size and minimum Lubell value full squashed flat antichain is maximal, and it is as well the unique volume minimizer.
\item[$n=2k-2$] There is a unique full squashed flat antichain with minimum Lubell value, namely, the full squashed flat antichain with minimum volume, which is also maximal.
    It is the first maximal (and full) squashed flat antichain of minimum size.
\item[$n\leq 2k-3$] All full squashed flat antichains with minimum Lubell value are maximal, and each occurs before the first minimum size full squashed flat antichain.
\end{description}

Let $L(n,k)$ be the minimum Lubell function value over all maximal squashed flat antichains
in $\binom{[n]}k\cup\binom{[n]}{k-1}$.
As the optimal full and maximal squashed flat antichains differ only marginally,
it is easy to verify that the following asymptotic result still holds for full squashed flat antichains.
For brevity, we only look at maximal squashed flat antichains here.

\begin{corollary}
For fixed $k\ge 1$ one has\/ $\displaystyle \lim_{n\to\infty}L(n,k)=1-\frac{(k-1)^{k-1}}{k^k}$.
\end{corollary}

\begin{proof}
For the limit to be shown, we can assume that $\lambda=(n-k+1)/k$ is large.
Considering this, Corollary \ref{minMFSA} implies that for an optimal
maximal squashed flat antichain we can choose
$a_1=0$ and $a_i=\lfloor(i-1)(\lambda+1)\rfloor$ for $i=2,\dots,k$.

For $2\le i\le k-1$ we have
$ \textstyle a_i~=~\lfloor (i-1)(\lambda+1)\rfloor ~=~ \lfloor \frac{i-1}k (n+1)\rfloor ~\le~ \frac ik n < n$.
Consequently, for $i\ne k$ the summands
\[ \alpha\binom{a_i}i-\beta\binom{a_i}{i-1}=\frac{\binom{a_i}i}{\binom nk}-\frac{\binom{a_i}{i-1}}{\binom n{k-1}} \]
on the right-hand side of (\ref{minweight}) all tend to $0$ as $n\to\infty$.

The claim follows by $\beta\binom n{k-1}=1$ and the by fact that
\begin{eqnarray*}
\alpha\binom{a_k}k-\beta\binom{a_k}{k-1} & = &
\frac{\binom{\lfloor\frac{k-1}k(n+1)\rfloor}k}{\binom nk}
-\frac{\binom{\lfloor\frac{k-1}k(n+1)\rfloor}{k-1}}{\binom n{k-1}} \\[1ex]
& = & \left(\frac{\lfloor\frac{k-1}k(n+1)\rfloor-k+1}{n-k+1}-1\right)
\prod_{j=0}^{k-2}\frac{\lfloor\frac{k-1}k(n+1)\rfloor-j}{n-j},
\end{eqnarray*}
which as $n\to\infty$ tends to
\[ \left(\frac{k-1}k-1\right)\left(\frac{k-1}k\right)^{k-1}~=~-\frac{(k-1)^{k-1}}{k^k}.\]
\end{proof}
%%%%%%%%%%%%%%%%%%%%%%%%%%%%%%%%%%%%%%%%%%%%%%%%%%%ZZZZZZZZZZZZZZ  end section cases of special interest
\section{Open Problems}

Here are a couple of new problems to expand on the results in this paper.

\begin{description}
%old version
% 1. What are the least (greatest) values of $m$ to achieve minimum size, minimum volume, and minimum BLYM value, respectively?
% Characterise the last $k$-set in each case.
% 2. How many different MSFA achieve each of minimum size, minimum volume, and minimum BLYM value?
% When are such MSFA unique?
% 3. When does an MSFA simultaneously achieve more than one of minimum size, minimum volume, and minimum BLYM value?
%\item[Problem 1.]
%Characterise those maximal squashed flat antichains that simultaneously attain
%minimum size and minimum Lubell function value.
% 4. It is well known that the largest antichain(s) in the Boolean Lattice are the maximum uniform antichain(s) in the
% middle level(s). Determine the maximum size MSFA for each $(n,k)$ and the size values that can be achieved
% between the  minimum and maximum size MSFA.
% 5. Which, if any MSFA, has an equivalent antichain which is not flat?
\item[Problem 1.]
Is there a maximal squashed flat antichain which has an equivalent
non-flat antichain (w.r.t.~the equivalence relation in the introduction)?
We conjecture the answer to be negative.
In particular, and noting the result in Section~4.2, is there a pair of values of $n,k$ such that the unique maximal squashed flat antichain with minimum size and minimum volume has an equivalent non-flat antichain?

\item[Problem 2.]
Extend the results here for full and maximal squashed flat antichains to antichains on two levels $k$ and $l$ that are
not necessarily consecutive.
%A generalised maximal squashed flat antichain (GMSFA) is a maximal squashed antichain which contains sets of exactly two sizes which do
%not need to be consecutive. Determine the minimum size, volume, and Lubell function values and answer the above
%questions in this more general setting.
%\item[Problem 4.]
%Consider the analogous questions for FSFA.
\end{description}

%%%%%%%%%%%%%%%%%%%%%%%%%%%%%%%%%%%%%%%%%%%%%%%%%%%ZZZZZZZZZZZZZZ  required remark
\section{Data Sharing}

Data sharing not applicable to this article as no datasets were generated or analysed during the current study.

%%%%%%%%%%%%%%%%%%%%%%%%%%%%%%%%%%%%%%%%%%%%%%%%%%%ZZZZZZZZZZZZZZ  end section oprn problems
%%\section{Acknowledgements}\label{sec:Concl}
%%%%%%%%%%%%%%%%%%%%%%%%%%%%%%%%%%%%%%%%%%%%%%%%%%%%%%%%%%%
%   we thank the co-authors for their infinite patience through
%   the delay in preparing this revision
%%%%%%%%%%%%%%%%%%%%%%%%%%%%%%%%%%%%%%%%%%%%

%Research supported in part by a grant from the Simons Foundation
%(\#282896 to Jerrold Griggs).

%%%%%%%%%%%%%%%%%%%%%%%%%%%%%%%%%%%%%%%%%%%%%%%%%%%%%%%%%%%%%%%%%%%%%%%%%%%%%%%%%%%%%%%%%%%%%
%%%%%%%%%%%%%%%%%%%%%%%%%%%%%%%%%%%%%%%%%%%%%%%%%%%%%%%%%%%%%%%%%%%%%%%%%%%%%%%%%%%%%%%%%%%%%

\end{document}